\documentclass[11pt]{amsart} 
\usepackage{amssymb, amscd}
\usepackage{latexsym,epsfig}
\usepackage{pst-all}
\numberwithin{equation}{section}

\newcommand{\VVl}{\VV_{\lambda}} 
\newcommand{\Pa}{\lambda = \{\lambda_1 \geq \lambda_2 \geq \lambda_3\geq 0\}} 
\newcommand{\lam}{\lambda}
\newcommand\la[1]{\lambda_{{#1}}} 
\newcommand{\Sy}{{\mathrm{Sym}}} \newcommand{\we}{\wedge}
\newcommand{\noir}{\Sy^{\la1-\la2}\VV\otimes
\Sy^{\la2-\la3}\left(\we^2\VV\right) \otimes
\Sy^{\la3}\left(\we^3\VV\right)} 
\newcommand{\VV}{\mathbb V}

\newcommand{\FF}{\mathbb F}
\newcommand{\eul}{e_c(\hy,\VVl)} 
\newcommand\cohi{\dim H_c^{i}(\hy,\VVl)}

\newcommand{\ls}{R^1\pi_*(\mathbb Q)} \newcommand{\hy}{{\mathcal H}_3}
\newcommand{\md}{{\mathcal M}_3} \newcommand{\mde}{{\mathcal M}_{3,1}}

\newcommand{\Sc}{{\mathrm{SL}}(2,{\mathbb C}) \times {\mathbb C}^*}
\newcommand{\es}{{\mathrm{SL}}(2, {\mathbb C})} \newcommand{\CC}{\mathbb C}
\newcommand{\e}{\varepsilon} 
\newcommand{\QQ}{\mathbb Q}
\newcommand{\einv}{\varepsilon^{-1}}
\newcommand{\Sp}{{\mathrm{Sp}}(6,\QQ)} 

\newcommand{\mat}{\left(\begin{array}{cc} a & b \\ c & d \end{array}\right)}

\newcommand{\uno}{\textbf{1}}
\newcommand{\Si}{\Sigma}
\newcommand{\pmenoe}{\left({\mathbb P}^1/V_4 - E\right)}

\newtheorem{theorem}{Theorem}[section]

\newtheorem{proposition}[theorem]{Proposition}

\newtheorem{definition-lemma}[theorem]{Definition-Lemma}

\theoremstyle{definition}

\theoremstyle{remark}

\begin{document}

\title[]{The euler characteristic of local systems \\ on the moduli
 of genus $3$ \\
hyperelliptic curves} \author{Gilberto Bini} 
\address{Dipartimento di Matematica, Universit\`a degli Studi di Milano,
Via C.\ Saldini, 50, 20133 Milano, Italia}
\email{Gilberto.Bini@mat.unimi.it}

\author{Gerard van der Geer}
\address{Faculteit Wiskunde en Informatica, University of
Amsterdam, Plantage Muidergracht 24, 1018 TV Amsterdam, The Netherlands.}
\email{geer@science.uva.nl}
%
%



\begin{abstract}
For a partition $\Pa$ of non-negative integers, we calculate the
Euler characteristic of the local system $\VVl$ on the
moduli space of genus $3$ hyperelliptic curves using a suitable stratification. 
For some $\lam$ of low degree, we
make a guess for the motivic Euler characteristic of $\VVl$ using
counting curves over finite fields.

\end{abstract}

\maketitle

\begin{section}{Introduction}
\label{sec: intro}

Let $\hy$ be the moduli space of genus $3$ hyperelliptic curves. It is
a $5$-dimensional substack of the Deligne-Mumford stack $\md$ of
smooth curves of genus $3$. 
The universal curve $\pi: \mde \rightarrow \md$ defines 
a natural local system $\ls$ of rank $6$ on $\md$. It comes with
a non-degenerate symplectic pairing.
The inclusion morphism $\iota: \hy
\rightarrow \md$ defines a natural local system $\VV:=\iota^*(\ls)$ on $\hy$. 

For any partition $\Pa$ of weight $|\lambda|=\lambda_1+
\lambda_2+\lambda_3$, consider
the irreducible representation of $\Sp$ associated with $\lam$. Any
such representation yields a symplectic local system $\VVl$ on 
$\hy$,
which appears `for the first time' in the decomposition of
$$
\noir.
$$
If, for example, $\lambda= \{\la1 \geq 0\geq 0\}$, then
$\VV_{\lambda}=\Sy^{\la1}(\VV)$. 

The cohomology with compact support of $\hy$ with local coefficients
in $\VV_{\lambda}$ is supposed to give interesting motives related
to automorphic forms. As a first step in understanding this cohomology
one wants to know the Euler characteristic of $\VV_{\lambda}$. This
was calculated for genus $2$ by Getzler in \cite{get}.
 In the present paper we calculate the
Euler characteristic
$$
\label{eulcar}
\eul=\sum_{i=0}^{10}(-1)^i\cohi
$$
for any local system $\VVl$ on $\hy$. We do this by using a
stratification of $\hy\otimes \CC$ by a
union of quasi-projective varieties $\Sigma(G)$, where $G$ is a finite
subgroup of $\Sc$, which acts on $\VVl$. By standard properties of the
Euler characteristic of local systems, we thus have
$$
\eul = \sum_G e_c(\Sigma(G))\dim(\VVl^G),
$$ where $e_c(\Sigma(G))$ is the topological Euler characteristic of
$\Sigma(G)$ and $\VVl^G$ is the space of $G$-invariants. We determine
$e_c(\Sigma(G))$ via elementary topological arguments and
$\dim(\VVl^G)$ via character theory. Getzler wrote down the generating series
of Euler characteristics in \cite{get}; however for genus $2$ already
this leads to unwieldy
rational functions. We give a short algorithm 
that calculates these number efficiently. 

This calculation is a step in the program
to understand the motivic Euler characteristic
$$
\sum_{i=0}^{10}(-1)^i[H^i_c(\hy,\VVl)],
$$
where $[H^i_c(\hy,\VVl)]$ is the class of the cohomology with
compact support in the Grothendieck ring of mixed $\QQ$-Hodge
structures. The hope is that in analogy to the genus $2$ case 
(cf.\ \cite{vand}), 
one could use this motivic Euler characteristic to describe properties of 
Siegel modular forms of genus $3$, of which very little is known. 
In Section \ref{sec4}, we provide some conjectural formulas of 
the motivic Euler characteristic
for specific low values of $|\lam|$ based on calculations
over finite fields.

Throughout the paper, $\e_n$ denotes a primitive $n$-th root of
unity.

\end{section}

\begin{section}{Stabilizers of hyperelliptic curves}
\label{sec1}

Let $C$ be a hyperelliptic curve of genus $3$ over the field of
complex numbers $\CC$.  
Then $C$ is a degree two cover of ${\mathbb P}^1$ with eight ramification points.
It can be given as a curve in the
$(X,Y)$-plane by an equation of the form $Y^2=f(X)$, where
$f(X)$ is a polynomial in $\CC[X]$ of degree $7$ or $8$.

The group $\Sc$ acts on the $(X,Y)$-plane as
follows. An element
$$
(A, \xi)= \left(\mat, \xi \right) \in \Sc
$$
acts via
$$
\label{actie1}
(A,\xi)\cdot (X,Y):= \left(\frac{aX+b}{cX+d}, \frac{\xi
Y}{(cX+d)^4}\right).
$$
Suppose that $G \leq\Sc$ stabilizes $C$. 
Consider the image $G'$ of $G$ under the projection of $\Sc$ onto 
$SL(2, {\mathbb C})$.
Clearly, $G'$ acts as
a group of rational transformations on the complex projective line. It also
permutes the set of ramification points of $C$. Note that the kernel
of this action is the subgroup generated by the central element
$-I$. 
By the classification of finite subgroups of $\es$ (see \cite{kl}), $G'$ 
must be isomorphic to one of the following groups:
\begin{itemize}
\item [i)]\label{cyclic} the cyclic group $C_n$ of order $n=2, 4, 14$;
\item [ii)] \label{quat} the quaternionic group $Q_{4n}$ of order
$4n=8, 12, 16, 24, 32$;
\item [iii)]\label{tetr} the group $O$ of symmetries of a cube.
\end{itemize}

For the purposes of what follows, we briefly recall the presentation
of the groups in i), ii), iii) as subgroups of $\es$. Any cyclic group
of order $n$ in $\es$ is conjugated to the group generated by the
matrix
$$
\left(
\begin{array}{cc}
\e_n & 0 \\
0 & \einv_n
\end{array}
\right).
$$ Any quaternionic subgroup of order $4n$, $n \geq 2$, is conjugated
to the group with generators
$$
S=
\left(
\begin{array}{cc}
\e_{2n} & 0 \\
0 & \einv_{2n}
\end{array}
\right) \quad \hbox{\rm and} \quad
U=
\left(
\begin{array}{cc}
0 & 1 \\
-1 & 0
\end{array}
\right). 
$$ 
Finally, the group $O$ is conjugated to the group generated by the matrices
$$
T=
\frac{-1}{\sqrt{2}}\left(
\begin{array}{cc}
1 & \e_8 \\
\e_8^3 & 1
\end{array}
\right) \quad \hbox{\rm and} \quad
U=
\left(
\begin{array}{cc}
0 & 1 \\
-1 & 0
\end{array}
\right). 
$$
Remarkably, the isomorphism type of $G'$ determines the whole
structure of $G$. Indeed, for any matrix $A \in G'$ there exist
two non-zero complex numbers $\pm \xi$ such that
\begin{equation}
\label{character}
\xi^2 Y^2=(cX+d)^8f\left(\frac{aX+b}{CX+d}\right),
\end{equation}
where
$$ 
A= \mat.
$$
The assignment
$$
u: G' \rightarrow  {\mathbb C}^*, \quad A \mapsto  \xi^2,
$$ 
is a character of a one-dimensional representation of $G'$
because $u(I)=1$. Thus,
the group $G \leq \Sc$ contains all pairs $(A,\pm u(A))$, where $A$
varies in one of the groups $G'$ listed in i), ii), iii), and $u$ is a
one-dimensional character of $G'$ that satisfies
\eqref{character}. Hence, $\#G= 2\#G'$. 

As a consequence, there are only finitely many non-isomorphic groups
$G$ which arise as possible stabilizers of genus $3$ hyperelliptic
curves. Each of them induces a permutation action on a set of eight
points in ${\mathbb P}^1$. Thus, we can deduce a \emph{normal form} of
curves which are stabilized by $G$. Examples and explicit
computations can be found, for instance, in \cite{sha}.
There, the stabilizers 
are not described as subgroups of $PSL(2, {\mathbb C}) \times {\mathbb C}^*$. 
 It is however easy to verify a correspondence between the two descriptions.

In Table 1 we list all possible groups in terms of $G'$ and $u$, as
well as the associated normal form. To this end, we need to review
some conventional notation from character theory. In general, we shall
denote by ${\uno}$ the trivial character of $G'$. If $G'$ is the
cyclic group of order $n$, there are $n-1$ nontrivial characters
$\chi^k$ such that
$$
\chi^k \left( \left(
\begin{array}{cc}
\e_n & 0 \\
0 & \einv_n
\end{array}
\right)\right) )=\e^k_n, \quad 1 \leq k \leq n-1.
$$
On the other hand, the quaternionic group $Q_{4n}$ has only three
non-trivial characters of one-dimensional representations, namely:
\begin{displaymath}
\begin{array}{|c|r|r|}
\hline
\chi & \chi(S) & \chi(U)
  \\  \hline
\chi_0 & 1 & -1  \\   
\chi_+ & -1 & -i^n  \\  
\chi_- & -1 & i^n  \\ \hline 
\end{array}
\end{displaymath}
The group $O$ has a unique $1$-dimensional character $\rho$, which is
not trivial.
\medskip
\smallskip
\medskip\noindent
\vbox{
\bigskip\centerline{\def\quad{\hskip 0.6em\relax}
\def\quod{\hskip 0.5em\relax }
\vbox{\offinterlineskip
\hrule
\halign{&\vrule#&\strut\quod\hfil#\quad\cr
height2pt&\omit&&\omit&&\omit&\cr
&{\rm name} &&$(G',u)$&&{\rm Normal Form $Y^2=f(X)$}&\cr
height2pt&\omit&&\omit&&\omit&\cr
\noalign{\hrule}
height2pt&\omit&&\omit&&\omit&\cr
&$G_1$&& $(C_2, {\uno})$ &&$(X^2-1)(X^6+\sum_{i=1}^5a_iX^{6-i}+1)$&\cr
height2pt&\omit&&\omit&&\omit&\cr
&$G_2$&& $(C_4, {\uno})$ &&$X^8+b_1X^6+b_2X^4+b_3X^2+1$&\cr
height2pt&\omit&&\omit&&\omit&\cr
&$G_3$&& $(Q_8, {\uno})$ &&$(X^4+c_1X^2+1)(X^8+c_2X^4+1)$&\cr
height2pt&\omit&&\omit&&\omit&\cr
&$G_4$&& $(C_4, {\chi^2})$ &&$X(X^6+d_1X^4+d_2X^2+1)$&\cr
height2pt&\omit&&\omit&&\omit&\cr
&$G_5$&& $(Q_{16}, {\uno})$ &&$X^8+fX^4+1$&\cr
height2pt&\omit&&\omit&&\omit&\cr
&$G_6$&& $(Q_8, {\chi_0})$ &&$(X^4-1)(X^4+l X^2+1)$&\cr
height2pt&\omit&&\omit&&\omit&\cr
&$G_7$&& $(Q_{12}, {\uno})$ &&$X(X^6+mX^3+1)$&\cr
height2pt&\omit&&\omit&&\omit&\cr
&$G_8$&& $(Q_{32}, {\chi_{-}})$ &&$X^8-1$&\cr
height2pt&\omit&&\omit&&\omit&\cr
&$G_9$&& $(O, {\uno})$ &&$X^8+14\, X^4+1$&\cr
height2pt&\omit&&\omit&&\omit&\cr
&$G_{10}$&& $(Q_{24}, {\chi_{-}})$ &&$X(X^6-1)$&\cr
height2pt&\omit&&\omit&&\omit&\cr
&$G_{11}$&& $(C_{14}, {\chi^6})$ &&$X^7-1$&\cr
height2pt&\omit&&\omit&&\omit&\cr
} \hrule}
}}
\centerline{\sc Table 1: \rm Groups and Associated Normal Forms}
\bigskip

We remark that the normal forms in Table 1 are equivalent to the
equations given in \cite{sha}, Table 3. For example, the
map
$$ (X,Y) \mapsto \left(\frac{-iX+i}{X+1},
\frac{\sqrt{8}\, \e_8}{\sqrt{2-l}(X+1)^4}\right), \quad l\neq 2,
$$ 
transforms the normal form associated with $G_6$ to
$$
\label{altnr}
Y^2=X(X^2-1)(X^4+LX^2+1),
$$
where $L= -(12+2l)/(2-l)$.
Additionally, the character $u$ changes too. However, this does not affect
the calculation of $\eul$ - see Section~\ref{sec3}.
\end{section}

\begin{section}{The stratification of $\hy$}
\label{sec2}

For each group $G_i$ in Table 1, define $\Sigma_i$ to be the locally
closed sublocus of $\hy$ that contains all curves $C$ whose stabilizer
is \emph{exactly} $G_i$. As seen in Section \ref{sec1}, the
corresponding group $G'_i$ is a permutation group on a set of eight
elements. We thus obtain a stratification of $\hy$ if the relation
$G_i' \leq G_j'$ is interpreted as an inclusion of permutation
groups. In other words, $G'_i$ is a subgroup of $G_j'$, and any set
of eight elements, which is permuted by $G'_i$, can be decomposed in
$G'_j$-orbits. All possible relations are displayed 
in the diagram below.
\par
\bigskip

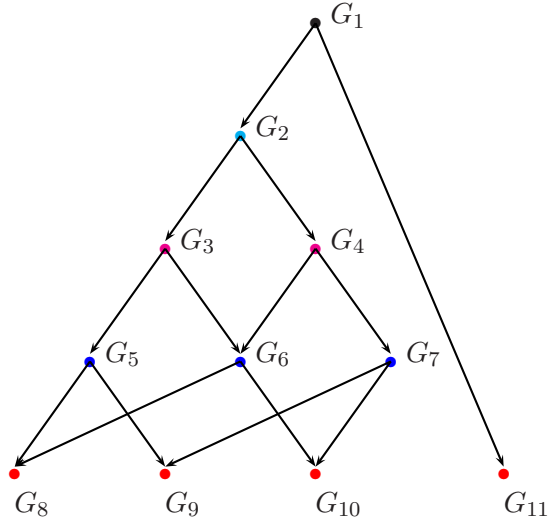
\begin{figure}
\begin{center}
\setlength{\unitlength}{0.1cm}
\begin{picture}(10,10)(40,80)
\put(10,-5){$G_8$}
\put(30,-5){$G_9$}
\put(50,-5){$G_{10}$}
\put(75,-5){$G_{11}$}
\put(22,15){$G_5$}
\put(42,15){$G_6$}
\put(62,15){$G_7$}
\put(32,30){$G_3$}
\put(52,30){$G_4$}
\put(42,45){$G_2$}
\put(52,60){$G_1$}
\put(10,0){\color{red}\circle*{1.5}}
\put(30,0){\color{red}\circle*{1.5}}
\put(50,0){\color{red}\circle*{1.5}}
\put(75,0){\color{red}\circle*{1.5}}
\put(20,15){\color{blue}\circle*{1.5}}
\put(40,15){\color{blue}\circle*{1.5}}
\put(60,15){\color{blue}\circle*{1.5}}
\put(30,30){\color{magenta}\circle*{1.5}}
\put(50,30){\color{magenta}\circle*{1.5}}
\put(40,45){\color{cyan}\circle*{1.5}}
\put(50,60){\circle*{1.5}}
\psline{->}(2,1.5)(1,0.1)
\psline{->}(2,1.5)(3,0.1)
\psline{->}(4,1.5)(1,0.1)
\psline{->}(4,1.5)(5,0.1)
\psline{->}(6,1.5)(3,0.1)
\psline{->}(6,1.5)(5,0.1)
\psline{->}(5,6)(7.5,0.1)
\psline{->}(3,3)(2,1.6)
\psline{->}(3,3)(4,1.6)
\psline{->}(5,3)(4,1.6)
\psline{->}(5,3)(6,1.6)
\psline{->}(4,4.5)(3,3.1)
\psline{->}(4,4.5)(5,3.1)
\psline{->}(5,6)(4,4.6)
\end{picture}
\end{center}
\end{figure}

\vskip 6cm
\par

\vskip 1. truecm
From this diagram, to be justified later,
we also deduce information on the strata
$\Sigma_i$. Actually, we have:
\begin{enumerate}
\item  $\hy=\bigcup_{i=1}^{11}\Si_i$;
\item  $\Si_i \cap \Si_j = \emptyset$ for $i \neq j$;
\item $\Si_j \subseteq \overline{\Si}_i$ whenever $G'_i \leq G'_j$.
\end{enumerate}

As explained in Section \ref{sec: intro}, we need to calculate the
topological Euler characteristic $e_c$ of all the strata. Since
$e_c(\hy)=1$, we work out $e(\Si_i)$, $i=2,\ldots, 10$ and we deduce
$e_c(\Sigma_1)$.

\textbf{0-dimensional strata.} The stratum $\Sigma_i$ for $i=8,9,10,11$
is clearly $0$-dimensional and irreducible, so its Euler number is $1$.

\textbf{1-dimensional strata.} The strata corresponding to
$G_5$, $G_6$, $G_7$ are $1$-dimensional. Moreover, let us consider the
following subsets of ${\mathbb P}^1$:
$${\mathcal O}_1:=\{\e_8^k; \, 0 \leq k \leq 7\},
$$
$$
{\mathcal O}_2:=\{0, \infty, \pm 1, \pm \e_3, \pm \e_3^2\},
$$
$$ 
{\mathcal O}_3:=\{\pm \alpha_1, \pm i\alpha_1, \pm 1/\alpha_1, \pm
i/\alpha_1\},
$$ where $\alpha_1$ is a root of the polynomial
$X^2-(i+1)X-i$.

It is easy to verify that ${\mathcal O}_1$ is a $G_8'$-orbit, a union
of two $G'_5$-orbits and a union of three $G'_6$-orbits. On the other
hand, ${\mathcal O}_2$ is a union of three $G'_7$-orbits, a union of
three $G'_6$-orbits and a union of two $G'_{10}$-orbits.  Finally,
${\mathcal O}_3$ is a full $G'_{9}$-orbit, a union of two
$G'_5$-orbits and a union of three $G'_7$-orbits. This justifies the
lower row of directed edges in the above diagram.

As for the Euler number $e_c$, the following holds.
\begin{proposition}
\label{onedim}
The topological Euler characteristic of $\Sigma_i$, $i=5,6,7$, is
equal to $-2$.\end{proposition} 
\proof We just prove the statement for
$\Sigma_5$, the other cases being similar.  For $f \in \CC -\{\pm 2\}$,
consider the set of hyperelliptic curves $C_{f}$ with equation
$Y^2=X^8+fX^4+1$. By direct inspection, two such curves $C_{f_1}$ and
$C_{f_2}$ are isomorphic if and only if $f_1=\pm f_2$. Note that
$\Si_8$ and $\Si_9$ are the isomorphism classes of $C_{0}$ and
$C_{14}$, respectively. Therefore, there exists an isomorphism
$\Phi:\Sigma_{5}\cup\Si_8 \cup \Si_9 \rightarrow \CC - \{4\}$ which
maps the orbit of $C_f$ to $f^2$. Accordingly, the topological Euler
characteristic of $\Si_5$ is $-2$.  \qed

\textbf{2-dimensional strata.} As readily checked from Table 1, the
strata corresponding to $G_3$ and $G_4$ have dimension two. It is easy
to deduce from the ramification sets in  ${\mathbb P}^1$ 
 that the following holds:
$$
\Si_5 \subset \overline{\Si}_3, \quad \Si_6 \subset \overline{\Si}_3,
$$
$$
\Si_6 \subset \overline{\Si}_4, \quad \Si_7 \subset \overline{\Si}_4.
$$

On the other hand, note that $\Si_5$ does not lie in the closure of
$\Si_4$. Equivalently, there is no set ${\mathcal S}$ of eight
elements which is both a union of $G'_4$-orbits and
$G'_5$-orbits. Indeed, any set ${\mathcal S}\subset {\mathbb P}^1$
has always two orbits of length one under the action of
$G'_4$. Conversely, the permutation action of $G'_5$ does not have any
fixed point.
\begin{proposition}
\label{twodim1}
The topological Euler characteristic of $\Si_3$ is $1$.
\end{proposition}  
\proof 
The group $G_3$ corresponds to the pair $(G'_3,\uno)$, where $G'_3$ is the
quaternionic group $Q_4 \cong C_2 \times C_2$. The group $G_3'$ induces a
permutation action on ${\mathbb P}^1$ via the group $V_4$ generated by
the transformations $x \mapsto - x$ and $ x \mapsto 1/x$. Denote by
$V(x)$ the orbit of $x$ under $V_4$. Note $\# V(a)=4$ unless $a \in \{
0,\infty,1,-1,i,-i\}$.

We recall that the normal form associated with $G_3$ is
\begin{equation}
\label{verk0}
Y^2=f(X)=(X^4+c_1X^2+1)(X^4+c_2X^2+1).
\end{equation}
Moreover, we have
\begin{equation}
\label{insieme} 
\{x: f(x)=0\} =\{\pm q_1, \pm 1/q_1, \pm q_2, \pm 1/q_2\},
\end{equation}
for distinct $q_1, q_2$ such that $\#V(q_1)=\#V(q_2)=4$. Note that
$c_i=-q_i^2-1/q_i^2$ for $i=1,2$.

Let $\{Y^2=f_1(X)\}$ and $\{Y^2=f_2(X)\}$ be two curves with
stabilizer $G_3$. They are isomorphic if and only if there exists a
rational transformation that maps $\{z: f_1(z)=0\}$ to
$\{z: f_2(z)=0\}$. All such transformations commute
with the elements of $V_4$. Therefore, two curves are isomorphic if
and only if there exists an automorphism of ${\mathbb P}^1/V_4$ which
preserves the set $E:=\{V(0),V(1),V(i)\}$, i.e.\  the ramification set
of ${\mathbb P}^1 \to {\mathbb P}^1/V_4$. Observe that the map
${\mathbb P}^1 \rightarrow {\mathbb P}^1/V_4$ sends $y$ to
$(y^2+1/y^2)/2$.

A curve $C$ with equation \eqref{verk0} has a larger stabilizer than
$G_3$ if and only if there exists $M \in {\rm SL}(2,\CC)$ - not in $G'_3$ -
which induces a permutation of \eqref{insieme} and a permutation of
the set $\{ 0,\infty,1,-1,i,-i\}$. By direct inspection, there is only
one possible $M$, namely:
$$
M=\left(
\begin{array}{cc}
\e_8 & 0 \\
0 & \einv_8
\end{array}
\right).
$$

In this case, $M$ induces the automorphism $x \mapsto ix$ on ${\mathbb
P}^1$ and the automorphism $z \mapsto -z$ on ${\mathbb P}^1/V_4$. Its
fixed points on ${\mathbb P}^1/V_4$ are $V(0)$ and $V(\e_8)$. 

Now, it is possible to give an alternative description of $\Sigma_3$,
which contains all curves whose stabilizer is \emph{exactly}
$G_3$. Denote by $\Delta$ the diagonal in $\pmenoe \times
\pmenoe$. Define a group $W_4$ of tranformations of ${\mathbb P}^1/V_4
\times {\mathbb P}^1/V_4$ as follows: $W_4$ is generated by
$\tau$, which interchanges both factors and $\iota$, which
simultaneously multiplies both factors by $i$. Note that $W_4$ is
isomorphic to the Klein four group. Therefore, $\Sigma_3$ can be
parametrized as
$$
\left(\pmenoe \times \pmenoe - \Delta -Z\right)/W_4,
$$
where 
$$
Z:=\{(V(a),V(ia): a\in \pmenoe - V(\varepsilon_8) \}.
$$

For the Euler number we get:
$$
e_c(\Sigma_3)= \frac14((-1)\times (-1) -(-1)-(-2))=1.
$$
\qed

\begin{proposition}
\label{twodim2}
The topological Euler characteristic of $\Si_4$ is $1$.
\end{proposition}
\proof The group $G_4$ corresponds to the pair $(G'_4, \chi^2)$, where
$G'_4$ is cyclic of order $2$. Now $G'_4$ induces a permutation action on
${\mathbb P}^1$ via the transformation $x \mapsto -x$. Denote by
$\sigma(x)$ the orbit of $x$ under such transformation.

We recall that the normal form associated with $G_4$ is
\begin{equation}
\label{verk}
Y^2=f(X)=X(X^6+d_1X^4+d_2X^2+1).
\end{equation}

Moreover, we have
$$ \{\infty\}\cup \{z: f(z)=0\} = \{\infty, 0, \pm a, \pm b, \pm c\}
$$ for some distinct $a,b,c \in \CC^*$. Therefore, any equation of the
form \eqref{verk} corresponds to the $5$-point set $\{\sigma(0),
\sigma(\infty), \sigma(a), \sigma(b), \sigma(c)\}$ on the ${\mathbb
P}^1$ which parametrizes the orbits $\{\sigma(x): x\in {\mathbb
P}^1\}$.

Let $\{Y^2=f_1(X)\}$ and $\{Y^2=f_2(X)\}$ be two curves with
stabilizer $G_4$. They are isomorphic if and only if there exists a
rational transformation that maps $\{\infty\}\cup \{z: f_1(z)=0\}$ to
$\{\infty\}\cup \{z: f_2(z)=0\}$ and fixes $0$ and $\infty$. Such a
transformation commutes with $x\rightarrow -x$. Consequently,
$\{Y^2=f_1(X)\}$ and $\{Y^2=f_2(X)\}$ are isomorphic if and only if
the associated $5$-point sets are mapped one onto the other by a
rational transformation which preserves $\sigma(0)$ and
$\sigma(\infty)$ and permutes the other three points. In other words,
an isomorphism class of curves with stabilizer $G_4$ defines an
element in ${\mathcal M}_{0,5}/{\mathfrak S}_3$, where ${\mathcal
M}_{0,5}$ is the moduli space of rational $5$-pointed curves and
${\mathfrak S}_3$ is the symmetric group of degree three. Conversely,
any element in ${\mathcal M}_{0,5}/{\mathfrak S}_3$ determines an
equivalence class of curves with stabilizer $G_4$.

Note that elements in ${\mathcal M}_{0,5}/{\mathfrak S}_3$ can be
written as $(0, \infty, 1, \sigma(u), \sigma(v))$ for some distinct $u, v\in {\mathbb P}^1-\{0,\infty, \pm 1\}$. The corresponding curve in $\Sigma_4$ have a bigger stabilizer if
and only if $\sigma(u)\sigma(v)=1$. As a consequence, $\Sigma_4$ can
be identified with ${\mathcal M}_{0,5}/{\mathfrak S}_3 - Y$, where $Y$
is the image of
$$
X:= \{(0, \infty, 1, \sigma(u), 1/\sigma(u))\} \subset {\mathcal M}_{0,5}
$$
under the quotient map onto ${\mathcal M}_{0,5}/{\mathfrak S}_3$.
Thus, we have
$$
e_c(\Sigma_4)= e_c({\mathcal M}_{0,5}/{\mathfrak S}_3)-e_c(Y)=1-e_c(Y)
$$
and
$$
e_c(X)=6e_c(Y)-r.
$$
Note that $e_c(X)= 2-4= -2$ since $\sigma(u) \notin \{\sigma(0),
\sigma(i), \sigma(1), \sigma(\infty)\}$. Additionally, $r=2$ since the
quotient map onto ${\mathcal M_{0,5}}/{\mathfrak S}_3$ is ramified
over $X$ when $\sigma(u)$ is the orbit of a primitive third root of
unity. Hence, the statement is completely proved.  \qed

\textbf{3-dimensional strata.} There is only a $3$-dimensional
stratum, namely $\Si_2$. As readily checked, both $\Sigma_3$ and
$\Sigma_4$ lie in the closure of $\Si_2$. 
\begin{proposition}
\label{threedim}
The topological Euler characteristic of $\Si_2$ is $2$.
\end{proposition}
\proof The group $G_2$ corresponds to the pair $(G'_2, \uno)$, where
$G'_2$ is cyclic of order two. As in Proposition \ref{twodim2}, $G'_2$
induces a permutation action on ${\mathbb P}^1$ via the transformation
$x \mapsto -x$. Again, denote by $\sigma(x)$ the orbit of $x$ under
such transformation.

We recall that the normal form associated with $G_2$ is
\begin{equation}
\label{verk1}
Y^2=f(X)=X^8+b_1X^6+b_2X^4+b_3X^2+1.
\end{equation}

Moreover, we have
$$ \{z: f(z)=0\} = \{\pm p_1, \pm p_2, \pm p_3 \pm p_4\}
$$ for some distinct $p_1,p_2,p_3, p_4 \in \CC^*$. Therefore, any
equation of the form \eqref{verk1} corresponds to the $4$-point set
$\{\sigma(p_1), \sigma(p_2), \sigma(p_3), \sigma(p_4)\}$ on the
${\mathbb P}^1$ which parametrizes the orbits $\{\sigma(x): x\in
{\mathbb P}^1\}$. 

Let $\{Y^2=f_1(X)\}$ and $\{Y^2=f_2(X)\}$ be two curves with
stabilizer $G_2$. They are isomorphic if and only if there exists a
rational transformation that maps $\{z: f_1(z)=0\}$ to
$\{z: f_2(z)=0\}$. All such possible transformations
commute with $x\rightarrow -x$. Consequently, $\{Y^2=f_1(X)\}$ and
$\{Y^2=f_2(X)\}$ are isomorphic if and only if the associated
$4$-point sets are mapped one onto the other by a rational
transformation. In other words, equivalence of curves with equation
\eqref{verk1} corresponds to equivalence of $4$-tuples of points in
${\mathbb P}^1$ under the action of ${\rm SL}(2,\CC)$ and the symmetric group of degree
$4$. Thus, an isomorphism class of curves stabilized by $G_2$ defines
a point in ${\mathcal M}_{0,4}/{\mathfrak S}_4$, where ${\mathcal
M}_{0,4}$ is the moduli space of $4$-pointed rational curves and
${\mathfrak S}_4$ is the symmetric group of order $4$. Note that
$e_c({\mathcal M}_{0,4}/{\mathfrak S}_4)=1$: see, for instance,~\cite{bgp}.

We finally observe that $\Sigma_2$ is not the whole ${\mathcal
M}_{0,4}/{\mathfrak S}_4$. In fact, we need to disregard all curves
with extra automorphisms, i.e., the ones in lower dimensional
strata. Therefore.
\begin{eqnarray*}
e_c(\Si_2)&= &e_c({\mathcal M}_{0,4}/{\mathfrak S}_4)-
\sum_{i=3}^{10}e_c(\Si_i)\\ &=& 1- (-6+2+3)=2.
\end{eqnarray*}
\qed

In Table 2, we list the dimension and the topological Euler
characteristic of all the strata in $\hy$.
$$\label{ta2}
\begin{array}{|c|r|r|r|r|r|r|r|r|r|r|r|}
 \hline 
i & 1&2&3&4&5&6&7&8&9&10&11\cr
\hline
\dim(\Sigma_i) & 5&3&2&2&1&1&1&0&0&0&0\cr
e_c(\Sigma_i) & -1 & 2& 1&1&-2&-2&-2&1&1&1&1\cr
\hline
\end{array}
$$

\centerline{\sc Table 2: \rm Some Topological Invariants of the Strata $\Si_i$.}
\smallskip
\end{section}

\begin{section}{The calculation of $\eul$}
\label{sec3}

Let $\gamma_j:\Sigma_j \rightarrow \hy$ be the embedding of $\Sigma_j$ in $\hy$.
By the properties of the Euler characteristic of local systems, we have
$$
\eul=\sum_{j=1}^{11}e_c\left(\Si_j,\gamma_j^*\left(\VVl\right)\right).
\label{enst}
$$

On the other hand, $\gamma_j^*\left(\VVl\right)$ is a local system on
$\Si_j$ with respect to $G_j$. Hence, \eqref{enst} can be written as
$$
\eul=\sum_{j=1}^{11}e_c(\Si_j)\dim(\VVl^{G_j}),
$$
where $\VVl^{G_j}$ is the space of $G_j$-invariants. In Section
\ref{sec2}, we computed $e_c(\Si_j)$. Now, we work out the dimension
of the corresponding invariant subspaces.

By definition, the fibre of the local system $\VV_{(1,0,0)}$ over a
curve $C$ is given by the cohomology group $H^1(C; \QQ)$. $\VVl$ is
thus obtained from the ${\rm Sp}(6, \QQ)$-module $\VV_{(1,0,0)}$ by standard
construction in representation theory (cfr. \cite{fh}). Obviously, any
group $G$ in Table 1 acts on $\VV_{(1,0,0)}$. This action yields a
homomorphism $\eta:G \rightarrow {\rm Sp}(6, \QQ)$. Let $(A,\xi)$ be an
element in $G$, where $A$ is a matrix with eigenvalues $a$ and
$a^{-1}$. By Corollary 3 in \cite{get}, the eigenvalues of $\eta(g)$
are given by
$$
a^2\xi, \quad a^{-2}\xi^{-1}, \quad a^{-2}\xi, \quad  a^2\xi^{-1}, 
\quad  \xi, \quad  \xi^{-1}.
$$

As a consequence, it is possible to compute the dimension of the $G$-invariant subspace of $\VVl$ by elementary character theory. More specifically, let $J_d$ be the symmetric function
$$
J_d(x_1, x_2, x_3)=h_d(x_1,x_1^{-1},x_2,x_2^{-1},x_3,x_3^{-1}),
$$
where $h_d$ is the complete symmetric function in six variables. Moreover, for any $\{\lam= \lam_1 \geq \lam_2 \geq \lam_3 \geq 0\}$, we denote by $J_{\lam}$ the determinant of the $3 \times 3$ matrix whose $i$-th row is
$$
\begin{array}{ccc}
(J_{\lam_i-i+2} &  J_{\lam_i-i+2} +J_{\lam_i-i} & J_{\lam_i-i+3} +J_{\lam_i-i-1}).
\end{array}
$$

By Proposition 24.22 in \cite{fh}, the following holds:
$$
\dim\left(\VVl^{G}\right)= \frac{1}{\#G}\sum_{g\in G}J_{\lam}(a^2\xi, a^{-2}\xi, \xi).
$$
For each of the groups $G_i$ we can list the pairs $(a,\xi)$
that occur as $g$ runs through $G$. If $(a,\xi)$ occurs, then
$(a,-\xi)$, $(-a,\xi)$ and $(-a,-\xi)$ occur too. For each $G_i$ in
Table 3 we give 
a set $Y_i$ of cardinality $\# G_i/4$ of 
pairs $(a,\xi)$ with multiplicity (indicated by
an exponent). The set $Y_i$ has the following property.
If we replace $(a,\xi) \in Y_i$ by the $4$ elements $(\pm a, \pm \xi)$ 
we get all the pairs with multiplicity corresponding to the $g\in G$.
This is indicated by the notation $(\pm a, \pm \xi) \in Y_i$.
\begin{theorem}
The Euler characteristic $\eul$ is given by
$$
\eul=\sum_{i=1}^{11} \frac{e_c(\Sigma_i)}{\# G_i} \sum_{(\pm a,\pm \xi) \in  Y_r}
J_{\lam}(a^2\xi, a^{-2}\xi, \xi),
$$
where the Euler numbers $e(\Sigma_i)$ and the sets
$Y_i$ are given in Tables 2 and 3.
\end{theorem}
\smallskip
\begin{displaymath}
\begin{array}{|c|r|}\hline
Y_1 & (1,1) \\
Y_2 & (1,1), (i,1)\\
Y_3 & (1,1),(i,1)^3 \\
Y_4 & (1,1), (i,i) \\
Y_5 & (1,1), (\e_{16}^2,1),(\e_{16}^6,1),(i,1)^5\\
Y_6 & (1,1),(i,1),(i,i)^2 \\
Y_7 & (1,1),(\e_{12}^2,1),(\e_{12}^4,1),(i,1)^3 \\
Y_8 & (1,1),(\e_{16},i), (\e_{16}^2,1),(\e_{16}^3,i),
(\e_{16}^5,i), (\e_{16}^6,1), (\e_{16}^7,i), (i,i)^4,(i,1)^5\\
Y_{9} & (1,1), (i,1)^9, (\e_{12}^2,1)^4,(\e_{12}^4,1)^4,
(\e_{16}^2,1)^3, (\e_{16}^6,1)^3\\
Y_{10} & (1,1),(\e_{14},\e_{14}^3),(\e_{14}^2,\e_{14}^6),
(\e_{14}^3,\e_{14}^9),(\e_{14}^4,\e_{14}^{12}),
(\e_{14}^5,\e_{14}),(\e_{14}^6,\e_{14}^4)\\
Y_{11} & (1,1), (i,1)^9, (\e_{12}^2,1)^4,(\e_{12}^4,1)^4,
(\e_{16}^2,1)^3, (\e_{16}^6,1)^3\\
Y_{11} & (1,1),(\e_{12},i),(\e_{12}^5,i),(\e_{12}^2,1),
(\e_{12}^4,1),(i,i)^4,(i,1)^3\\
\hline
\end{array}
\end{displaymath}
\centerline{\hbox{\sc Table 3. \rm The Sets $Y_i$}}
\par
\smallskip

For example, the elements of the group $G_1$ are $(\pm {\rm Id},\pm 1)$.
If $\lam=(k,0,0)$, then the contribution from this group yields
\begin{eqnarray*}
\dim (\VV_{(k,0,0)}^{G_1})&=&
\frac{1}{4}\left\{2h_k(1,1,1,1,1,1)+
2h_k(-1,-1,-1,-1,-1,-1)\right\} \\
&=&\frac{1}{2}\binom{k+5}{k}\left(1+(-1)^k\right).
\end{eqnarray*}
In the following table we give the values of $\eul$ for all
$\lambda$ of weight $\leq 10$. Note that because of the hyperelliptic
involution $\eul =0$ if the weight is odd.
\smallskip
\begin{displaymath}
\begin{array}{|c|r||c|r|}\hline
(\lam_1, \lam_2, \lam_3) & \eul & (\lam_1, \lam_2, \lam_3) & \eul \\
\hline
(0,0,0) & 1  &(5,2,1) & -10 \\
(2,0,0) & -1 & (4,4,0) & -5\\
(1,1,0) & 0 & (4,3,1) & -4 \\
(4,0,0) & -1  & (4,2,2) & -7\\
(3,1,0) & 0 & (3,3,2) & -2\\
(2,2,0) & -1 &(10,0,0) & -17 \\
(2,1,1) & 0 & (9,1,0) & -22\\
(6,0,0) & -5 &(8,2,0) & -43 \\ 
(5,1,0) & -2 &(8,1,1) & -8\\
(4,2,0) & - 5 & (7,3,0) & -34\\
(4,1,1) & 0 &(7,2,1) & -32\\
(3,3,0) & 0 &(6,4,0) & -37\\
(3,2,1) & 0 &(6,3,1) & -26\\
(2,2,2) & -3 & (6,2,2) & -27\\
(8,0,0) & -7& (5,5,0)& -6\\
(7,1,0) & -8 & (5,4,1)& -22\\
(6,2,0) & -13 & (5,3,2)& -12\\
(6,1,1) & -2 & (4,4,2)& -15\\
(5,3,0) & -10 &(4,3,3) &  0\\
\hline
\end{array}
\end{displaymath}
\centerline{\sc Table 4: \rm Some Values of $e_c(\hy,\VVl)$}
\end{section}
\begin{section}{Some Remarks on the motivic Euler characteristic}
\label{sec4}
For partitions of small degree $|\lambda|$ it is not unreasonable
to expect that all cohomology of  $\VV_{\lambda}$ is Tate,
i.e., that the motivic Euler characteristic
$$
E_c(\hy,\VV_{\lambda})=
\sum_{i=0}^{10} (-1)^i [H_c^i(\hy,\VV_{\lambda})]
$$
is a polynomial in $L$, the Tate motive of weight $2$. It is well
known that $E_c(\hy,\VV_0)=L^5$.
One can calculate the trace of Frobenius on the $\ell$-adic
variant of $\VV_{\lambda}$ in characteristic $p$ on $\hy \otimes \FF_p$
by summing
$$
\sum_C {\rm Tr}(F,\VV_{\lambda}(H^1))/\# {\rm Aut}_{\FF_p}(C),
$$
where $C$ runs over a complete set of representatives of the
$\FF_p$-isomorphism classes of hyperelliptic curves of genus $3$
over $\FF_p$. We found that the following guesses for the motivic
Euler characteristic are compatible with these traces for $p=2,3$ and $5$
and with the values of $\eul$. 

\smallskip
\medskip\noindent
\vbox{
\bigskip\centerline{\def\quad{\hskip 0.6em\relax}
\def\quod{\hskip 0.5em\relax }
\vbox{\offinterlineskip
\hrule
\halign{&\vrule#&\strut\quod\hfil#\quad\cr
height2pt&\omit&&\omit&\cr
&$\lambda$ &&$E_c(\hy,\VV_{\lambda})$&\cr
height2pt&\omit&&\omit&\cr
\noalign{\hrule}
height2pt&\omit&&\omit&\cr
&$(0,0,0)$ && $L^5$& \cr
&$(2,0,0)$ && $-1$&\cr
&$(1,1,0)$ && $0$&\cr
&$(4,0,0)$ && $L^2-2$&\cr
&$(3,1,0)$ && $L^2-1$&\cr
&$(2,2,0)$ && $- L^6 +L^2-1$&\cr
&$(2,1,1)$ && $L^5 -L^4 -L^3 +L^2$&\cr
height2pt&\omit&&\omit&\cr
} \hrule}
}}
\centerline{\sc Table 5. \rm Motivic Euler Characteristics}
\bigskip

\end{section}

\end{document}